\documentclass[12pt]{amsart}
\usepackage{amsmath, amssymb, amsthm, bm, mathtools, enumitem, caption}
\usepackage{hyperref}

\usepackage[noabbrev,capitalize]{cleveref}
\usepackage{adjustbox}
\crefname{equation}{}{}
\usepackage{ytableau}
\usepackage{graphicx, tikz}
\usepackage[textheight=8.95in, textwidth=6.75in]{geometry}

\usepackage{microtype}

\newtheorem{theorem}{Theorem}[section]
\newtheorem{lemma}[theorem]{Lemma}
\newtheorem{corollary}[theorem]{Corollary}
\newtheorem{proposition}[theorem]{Proposition}

\newtheorem*{conjecture*}{Conjecture}

\theoremstyle{definition}

\newtheorem*{definition}{Definition}

\theoremstyle{remark}
\newtheorem*{remark}{Remark}

\newtheorem*{tworemarks}{Two Remarks}
\newtheorem*{example}{Example}

\numberwithin{equation}{section}

\newcommand{\R}{\mathbb R}

\DeclareMathOperator{\Tr}{Tr}

\newcommand{\C}{\mathbb C}
\newcommand{\SL}{\mathrm{SL}}

\newcommand{\Q}{\mathbb Q}

\newcommand{\Z}{\mathbb Z}

\title[Some topological genera and Jacobi forms]{Some topological genera and Jacobi forms}
\author{Tewodros Amdeberhan, Michael Griffin, \and Ken Ono}
\address{Dept. of Mathematics, Tulane University, New Orleans, LA 70118}
\email{tamdeber@tulane.edu}
\address{Dept. of Mathematics, Vanderbilt University, Nashville, TN 37240}
\email{michael.j.griffin@vanderbilt.edu}
\address{Dept. of Mathematics, University of Virginia, Charlottesville, VA 22904}
\email{ko5wk@virginia.edu}

\begin{document}

\begin{abstract} 
We revisit and elucidate the $\widehat{A}$-genus, Hirzebruch's $L$-genus and Witten's $W$-genus, cobordism invariants of special classes of manifolds. 
After slight modification, involving Hecke's trick, we find that the $\widehat{A}$-genus and $L$-genus arise directly from Jacobi's theta function.  For every  $k\geq 0,$ we obtain exact formulas for the quasimodular expressions of $\widehat{A}_k$ and $L_k$ as ``traces'' of partition Eisenstein series
$$
\widehat{\mathcal{A}}_k(\tau)= \Tr_k(\phi_{\widehat{A}};\tau)\ \ \ \ \ \ {\text {\rm and}}\ \ \ \ \ \
\mathcal{L}_k(\tau)= \Tr_k(\phi_L;\tau),
$$
which are easily converted to the original topological expressions.
Surprisingly, Ramanujan defined twists of the $\widehat{\mathcal{A}}_k(\tau)$ in his ``lost notebook'' in his study of derivatives of theta functions, decades before Borel and Hirzebruch rediscovered them in the context of spin manifolds. In addition, we show that the nonholomorphic $G_2^{\star}$-completion of the characteristic series of the Witten genus is the Jacobi theta function avatar of the $\widehat{A}$-genus.
\end{abstract}

\maketitle

\section{Introduction and Statement of Results}

A sequence of polynomials $f_1, f_2,\dots$ in the variables $p_1, p_2,\dots$ is {\it multiplicative} if the identity
$$
1+p_1 t+p_2t^2+\dots =(1+r_1t+r_2t^2+\dots)(1+s_1t+s_2t^2+\dots)
$$
implies that
$$
\sum_{n=1}^{\infty} f_n(p_1, p_2,\dots)t^n=\left( \sum_{a=1}^{\infty}f_a(r_1, r_2,\dots)t^a\right)
\left(\sum_{b=1}^{\infty}f_b(s_1, s_2,\dots)t^b\right).
$$
If $Q(z)$ is a power series with constant term 1, then one gets such sequences from the infinite product
\begin{equation}\label{infinite_product}
F(p_1, p_2,\dots;t):=\prod_{i=1}^{\infty} Q(x_i t)=1+f_1t+f_2t^2+\dots,
\end{equation}
where $p_k$ is the $k$th elementary symmetric function (in the variables, $x_1, x_2,\dots$) defined by
$$
p_k:=\sum_{i_1<i_2 <\dots < i_k} x_{i_1}x_{i_2}\cdots x_{i_k}.
$$

By work of Thom,  this combinatorial framework applies to the study of homomorphisms of cobordism rings of manifolds with prescribed structure.
The idea is that a {\it characteristic power series} $Q(z)$  encodes invariants of oriented manifolds, with dimensions that are multiples of 4, via its {\it genus} given by (\ref{infinite_product}).
Here the  $p_k$ represent the Pontryagin classes, the cohomology classes of real vector bundles.

We consider the number theoretic properties of some well-known examples (see these references for background \cite{Hirzebruch, HBJ}). We first consider the $\widehat{A}$-genus of spin manifolds discovered by Borel and Hirzebruch \cite{BorelHirzebruch, Hirzebruch}. The 
 first few values are
 \begin{displaymath}
\begin{split}
\widehat{A}_0&=1,\ \
\widehat{A}_1=-\frac{1}{24}p_1,\ \ \ 
\widehat{A}_2=\frac{1}{5760}\left(-4p_2+7p_1^2\right),\ \ \
\widehat{A}_3=\frac{1}{967680}\left(-16p_3+44p_1p_2-31p_1^3\right), \dots.
\end{split}
\end{displaymath}
This example is historically significant because of its role in the discovery of the Atiyah-Singer index theorem (for example, see Hitchin's expository article \cite{Hitchin}).
Atiyah and Singer discovered and employed their index theorem to explain the mysterious integrality of the values of the $\widehat{A}$-genus. To compute these values, they implemented (\ref{infinite_product}) with
\begin{equation} \label{Gen_sinc}
Q_{\widehat{A}}(z):=\frac{\frac{1}{2}\sqrt{z}}{\sinh(\frac{1}{2}\sqrt{z})}=1-\frac{z}{24}+\frac{7z^2}{5760}-\frac{31z^3}{967680}+\dots.
\end{equation}
Namely, the $\widehat{A}$ values (in order) are the coefficients of the formal power series 
\begin{eqnarray}
\widehat{A}(p_1, p_2,\dots; t)&=&\sum_{n=0}^{\infty} \widehat{A}_n t^{n}
=\prod_{i=1}^{\infty}  Q_{\widehat{A}}(x_i t)\nonumber\\
&=&1-\frac{1}{24}p_1t+\frac{1}{5760}(-4p_2+7p_1^2)t^2+\frac{1}{967680}(-16p_3+44p_1p_2-31p_1^3)t^3+\dots.\label{Aexample}
\end{eqnarray}

We prove that $\widehat{A}(p_1, p_2,\dots; t)$ (and the other genera in this paper), after minor modification, is essentially a  {\it Jacobi form} (see Chapter~2 of \cite{BFOR} or \cite{Zagier}) on $\C\times \mathbb{H}.$
This connection opens the door to new avenues of research, where the theory of modular forms can be brought to bear on the number theoretic properties of these topological invariants that arise as Fourier coefficients. Namely, it is natural to expect that this work will lead to new results about the asymptotic and congruence properties of topological genera.

To make this connection precise, we recall  the celebrated Jacobi theta function (see \cite{EichlerZagier, Zagier})
$$
\theta(z;\tau):=\sum_{n\in \Z}u^nq^{n^2/2},
$$
where $u:=e^{2\pi i z}$ and $q:=e^{2\pi i \tau}.$
This function is a Jacobi form for $\SL_2(\Z)$ of weight $1/2$ and index $1/2$.
We work instead with a slightly modified version of this function. Readers familiar with \cite{AGOS} should be aware that $\widetilde{\Theta}(z;\tau)$ here is slightly different from the one in that paper. Namely,  in terms of Dedekind's eta-function
$\eta(\tau):=q^{\frac1{24}}\prod_{n=1}^{\infty}(1-q^n),$
 it will be convenient for us to employ
\begin{equation}\label{ThetaDef}
    \widetilde{\Theta}(z;\tau)
     \ = \  \exp\left(\frac{\pi}{2}\cdot \frac{z^2}{\Im(\tau)}\right)\cdot u^{\frac{1}{2}}q^{\frac{1}{8}}\cdot \frac{\theta\left (z+\frac{\tau}{2}+\frac{1}{2}; \tau \right)}{\eta(\tau)^3}.
\end{equation}

We transform the $\widehat{A}$-genus, as described above, into the function
\begin{eqnarray}
\widehat{A}(X_\tau(s); t):=\prod_{x\in X_\tau(s)} Q_{\widehat{A}}(xt),\label{Ahat}
\end{eqnarray}
where  $\Im(\tau)>0$ and $s\in \R^{+},$ and
\begin{equation}
X_{\tau}(s):= \left \{\frac{1}{(m\tau+n)^2\cdot |m\tau +n|^s} \ : \ \gcd(m,n)=1\right \}.
\end{equation}
As a function on $\C\times \mathbb{H},$ we have the following identity in terms of the Jacobi theta function.

\begin{theorem}\label{Theorem1A}  We have that
$$
\lim_{s\rightarrow 0^+} \widehat{A}(X_{\tau}(s); (2\pi i z)^2 )= 2\pi i z\cdot \widetilde{\Theta}(z;\tau)^{-1}.
$$
\end{theorem}

The infinite product in (\ref{Ahat}) is taken over relatively prime pairs of  integers $(m,n)$ instead of $i=1, 2,\dots,$ as in (\ref{infinite_product}).  This modification does not lose any information. In fact, this reformulation will allow us to compute the values of the $\widehat{A}$-genus (and also $L$-genus)  as quasimodular forms (see Theorem~\ref{QuasimodularGenera}). Moreover, it is straightforward to reconstruct the original $\widehat{A}$-genus expressions (and also $L$-genus) using elementary properties of symmetric functions. The example after Theorem 1.4  illustrates this combinatorial and number theoretic procedure. Finally, we note that the dependence on $s> 0$ in the index set $X_{\tau}(s)$ is required, as we view these series as analytic functions, and  the introduction of $s$ guarantees convergence.

We also consider Hirzebruch's $L$-genus \cite{Hirzebruch}, which is the case of closed smooth oriented manifolds.  The first few values are
\begin{displaymath}
\begin{split}
{L}_0&=1,\ \
{L}_1=\frac{1}{3}p_1,\ \ \ 
{L}_2=\frac{1}{45}\left(7p_2-p_1^2\right),\ \ \
{L}_3=\frac{1}{945}\left(62p_3-13p_1p_2+2p_1^3\right), \dots.
\end{split}
\end{displaymath}
In terms of the characteristic power series
\begin{equation}\label{LGenus}
Q_{L}(z):=\frac{\sqrt{z}}{\tanh(\sqrt{z})}=1+\frac{z}{3}-\frac{z^2}{45}+\frac{2z^3}{945}-\dots,
\end{equation}
 the infinite product (\ref{infinite_product}) gives the generating function
 \begin{eqnarray}
L(p_1, p_2,\dots; t)&=&\sum_{n=0}^{\infty} L_n t^{n}  =\prod_{i=1}^{\infty}  Q_{L}(x_i t) \nonumber\\
&=&1+\frac{1}{3}p_1t+\frac{1}{45}(7p_2-p_1^2)t^2+\frac{1}{945}(62p_3-13p_1p_2+2p_1^3)t^3+\dots.\label{Lexample}
\end{eqnarray}

We prove that  $L(p_1, p_2,\dots; t)$, after minor modification, is also essentially a Jacobi form. As in the case of the 
$\widehat{A}$-genus, we transform the $L$-genus into the function
\begin{eqnarray}
L(X_\tau(s);t)&:=&\prod_{x\in X_{\tau}(s)} Q_L(xt). \label{Lhat}
\end{eqnarray}
As a function on $\C\times \mathbb{H},$ we have the following identity.

\begin{theorem}\label{Theorem1L} We have that
$$
\lim_{s\rightarrow 0^+} L(X_{\tau}(s); (\pi i z)^2) =\pi i z\cdot \frac{\widetilde{\Theta}(2z; \tau)}{\widetilde{\Theta}(z;\tau)^2}.
$$
\end{theorem}

Theorems~\ref{Theorem1A} and \ref{Theorem1L}  connect the $\widehat{A}$-genus and $L$-genus to the theory of  elliptic modular forms. As a corollary to Theorem~\ref{Theorem1A}, we relate the $\widehat{A}$-genus to the characteristic series of the {\it Witten genus} for compact oriented smooth spin manifolds with vanishing
first Pontryagin class, that naturally arises from modularity.
To make this precise, for
 integers $k\geq 1$ and $\Im(\tau)>0,$ the weight $2k$ Eisenstein series (see Ch. 1 of \cite{CBMS}) is
\begin{equation}
G_{2k}(\tau) := -\frac{B_{2k}}{2k}+2\sum_{n}\sigma_{2k-1}(n)q^n \ = \ \frac{(2k-1)!}{(2\pi i)^{2k}}\sum_{\substack{\omega \in \Z\oplus \Z\tau\\
\omega\neq 0}}\frac{1}{\omega^{2k}},
\end{equation}
where $B_{2k}$ is the  $2k$-th Bernoulli number and $\sigma_{\nu}(n):=\sum_{d\mid n}d^{\nu}.$ The first examples are
$$
G_2(\tau)=-\frac{1}{12}+2\sum_{n=1}^{\infty}\sigma_1(n)q^n,\ \ \
G_4(\tau)=\frac{1}{120}+2\sum_{n=1}^{\infty} \sigma_3(n)q^n,\ \ \
G_6(\tau)= -\frac{1}{252}+2\sum_{n=1}^{\infty}\sigma_5(n)q^n.
$$
Apart from $G_2,$ each $G_{2k}$ is a weight $2k$ holomorphic modular form on $\SL_2(\Z),$ and the
{\it quasimodular forms}  are the $q$-series in the polynomial ring
 (for example, see \cite{Zagier})
$$
\mathbb{C}[G_2,G_4,G_6]=\mathbb{C}[G_2, G_4, G_6, G_8, G_{10}, \dots].
$$

The modular Eisenstein series $G_4, G_6,\dots$ are compiled to form
the {\it Witten genus} \cite{Witten} (also see \cite{Hopkins}) via its characteristic series
\begin{equation}\label{WittenGenus}
Q_W(z)=\exp\left(\sum_{k\geq 2}\frac{G_{2k}(\tau)(2\pi i z)^{2k}}{(2k)!}\right).
\end{equation}
This identity implies that the Witten genus of a $4k$ dimensional compact oriented smooth spin manifold, with vanishing first Pontryagin class, is a weight $2k$ modular form with integral Fourier coefficients. It is natural to ask about the topological significance of the function that one obtains by including $G_2$ in this characteristic series.
It turns out that one obtains the Jacobi theta function avatar of the $\widehat{A}$-genus.

\begin{corollary}\label{WittenCorollary} We have that
$$
\lim_{s\rightarrow 0^+} \widehat{A}(X_{\tau}(s); (2\pi i z)^2)= \exp\left( (2\pi i  z)^2 \cdot \frac{G_2^{\star}(\tau)}{2}\right)\cdot Q_W(z),
$$
where $G_2^{\star}(\tau):=\frac{1}{4\pi\Im(\tau)}+G_2(\tau)$ is the nonholomorphic weight 2 modular Eisenstein series.
\end{corollary}

As a consequence of both Theorems~\ref{Theorem1A} and \ref{Theorem1L},  we obtain quasimodular representations of the $\widehat{A}$-genus and $L$-genus.  
We now turn to the problem of converting these quasimodular forms into the original topological expressions. To make this precise, we make the important observation that
these forms are given as 
traces of ``partition Eisenstein series,'' which are studied in \cite{AGOS, AOS, AOS2}.  
To define them, recall that a {\it partition of  a non-negative integer $k$} (see \cite{Andrews} for background on partitions) is any nonincreasing sequence of positive integers 
$$\lambda=(\lambda_1,\lambda_2,\dots, \lambda_s)$$ that sum to $k$, denoted $\lambda\vdash k.$ Equivalently, we let  $\lambda=(1^{m_1},\dots,k^{m_k})\vdash k$, where $m_j$ is the multiplicity of $j.$ Furthermore, the {\it length} of $\lambda$ is $\ell(\lambda):=m_1+\dots+m_k.$
For a partition $\lambda$, we
 define the weight $2k$ {\it partition Eisenstein series}
\begin{equation}\label{primary}
\lambda=(1^{m_1}, 2^{m_2},\dots, k^{m_k}) \vdash k \ \ \ \ \ \longmapsto \ \ \ \ \  G_{\lambda}(\tau):= G_2(\tau)^{m_1}G_4(\tau)^{m_2}\cdots G_{2k}(\tau)^{m_k}.
\end{equation}
In particular, the Eisenstein series $G_{2k}(\tau)$ corresponds to the partition $\lambda=(k).$ The $G_{\lambda}$ should not be mistaken for the partition Eisenstein series of Just and Schneider \cite{JustSchneider}.

If $\phi: \mathcal{P}\mapsto \C$ is a function on partitions, then for $k\geq 1$ we define the {\it partition Eisenstein trace}
\begin{equation}\label{PartitionTrace}
\Tr_k(\phi;\tau):=\sum_{\lambda \vdash k} \phi(\lambda)G_{\lambda}(\tau),
\end{equation}
which is a weight $2k$ quasimodular form.
By convention, for $k=0$, we let $\Tr_0(\phi;\tau):=1$.

We give quasimodular representations of  the $\widehat{A}$-genera and $L$-genera as partition Eisenstein traces.
To this end, we first note that $\widehat{A}(p_1, p_2,\dots; t)$ and $L(p_1, p_2,\dots; t)$ are of the form
$$
F(p_1, p_2,\dots; t)=1+\sum_{k=1}^{\infty} b_k(F; p_1,p_2,\dots)t^k,
$$ 
where each $b_k(F;p_1, p_2,\dots)$ is a homogeneous polynomial of {\it weighted} degree $k$. In other words, each monomial $p_1^{m_1}p_2^{m_2}\dots p_k^{m_k}$ has weight
$k=m_1+2m_2+\dots +km_k.$ This provides the unique representation
$$
b_k(F; p_1, p_2,\dots)= \widetilde{b}_k(F; s_1, s_2,\dots),
$$
where the $s_j:=x_1^j+x_2^j+\dots,$  are the $j$th power sum symmetric functions. Clearly, as a polynomial in $s_1, s_2,\dots, s_k,$ we have that
$\widetilde{b}_k(F;s_1, s_2,\dots)$ is also homogeneous of weighted degree $k.$ We simplify notation by associating partitions with monomials, where
\begin{equation}
  s_{\lambda}:=s_1^{m_1}s_2^{m_2}\cdots s_k^{m_k},
\end{equation}
 with $\lambda=(1^{m_1}, 2^{m_2},\dots, k^{m_k}) \vdash k$. Therefore, we have a decomposition
 $$
 b_k(F; p_1, p_2,\dots)=\widetilde{b}_k(F;s_1, s_2,\dots)=\sum_{\lambda \vdash k} \beta_F(\lambda)\cdot s_{\lambda}.
 $$
To each $b_k(F; p_1, p_2,\dots),$ we associate the weight $2k$ partition Eisenstein trace
\begin{equation}
\mathcal{F}_k(\tau):=\sum_{\lambda \vdash k} \beta^{\star}_F(\lambda)\cdot G_{\lambda}(\tau),
\end{equation}
where we modify the coefficients $\beta_F(\lambda)$ with a Bernoulli product as follows
\begin{equation}\label{BernoulliWeight}
\beta_F^{\star}(\lambda):= \beta_F(\lambda) \cdot \prod_{j=1}^{k}\left(\frac{2j}{B_{2j}}\right)^{m_j}.
\end{equation}

By letting $F=\widehat{A}(p_1,p_2,\dots; t)$ (resp. $F=L(p_1,p_2,\dots;t)),$  we obtain $\mathcal{\widehat{A}}_k(\tau)$ (resp. $\mathcal{L}_k(\tau)$),
the weight $2k$ quasimodular avatars of $\widehat{A}_k$ (resp. $L_k$).  To make this explicit, we define the functions
\begin{eqnarray}
\phi_{\widehat{A}}(\lambda)&:=&\prod_{j=1}^k\frac1{m_j!}\left(\frac{-1}{(2j)!}\right)^{m_j},  \label{phi_A}  \\
\phi_{L}(\lambda)&:=&\prod_{j=1}^k\frac1{m_j!}\left(\frac{4^j(4^j-2)}{(2j)!}\right)^{m_k}.
\end{eqnarray}
The following theorem gives the exact quasimodular expressions for these genera.

\begin{theorem}\label{QuasimodularGenera}
If $k$ is a positive integer, then as Fourier series we have
\begin{displaymath}
\begin{split}\widehat{\mathcal{A}}_k(\tau)= \Tr_k(\phi_{\widehat{A}};\tau),\\
\mathcal{L}_k(\tau)= \Tr_k(\phi_L;\tau).
\end{split}
\end{displaymath}
\end{theorem}

\begin{example} It is straightforward to derive the $\widehat{A}_k$ and $L_k$ (see (\ref{Aexample}) and (\ref{Lexample}))
using Theorem~\ref{QuasimodularGenera}. One transforms the quasimodular traces
$\Tr_k(\phi_{\widehat{A}};\tau)$ and $\Tr_k(\phi_L;\tau)$ into expressions in the power sum symmetric functions, and then, in turn, into expressions in the elementary symmetric functions. In view of (\ref{BernoulliWeight}), in the first step one
replaces each $G_{2j}(\tau)$ with $B_{2j}s_j/2j.$ 

For  $\widehat{A}_3$ and  $\widehat{A}_4,$ Theorem~\ref{QuasimodularGenera} gives
\begin{align*}
\widehat{\mathcal{A}}_3(\tau)&=
\Tr_3(\phi_{\widehat{A}};\tau)= \frac1{6!}(-G_6+15G_2G_4-15G_2^3),   \\
\widehat{\mathcal{A}}_4(\tau)&=
\Tr_4(\phi_{\widehat{A}};\tau)= \frac1{8!}(-G_8+28G_2G_6+35G_4^2-210G_2^2G_4+105G_2^4).
\end{align*}
After making the substitutions  $G_{2j}\mapsto B_{2j}s_j/2j$,  we apply the Newton-Gerard identities
$$
s_1=p_1, \ \  s_2=p_1^2-2p_2, \ \ s_3=p_1^3-3p_1p_2+3p_3, \ \  s_4=p_1^4-4p_1^2p_2+4p_1p_3+2p_2^2-4p_4,
$$
and we obtain
\begin{align*}
\widehat{A}_3(p_1,p_2,\dots)&=\frac1{967680}(-16p_3 +44p_1p_2 -31p_1^3), \\
\widehat{A}_4(p_1,p_2,\dots)&=\frac1{464486400}(-192p_4+512p_1p_3+208p_2^2 - 904p_1^2p_2 +381p_1^4).
\end{align*}
\end{example}

To our surprise, it turns out that
 Ramanujan discovered the quasimodular representations of the $\widehat{A}$-genus 100 years ago, decades before Borel and Hirzebruch resdiscovered them in the context of spin manifolds.
 In his ``lost notebook'', Ramanujan defined the $q$-series \cite[p. 369]{Rama}
\begin{equation}\label{U}
U_{2k}(q)=\frac{1^{2k+1}-3^{2k+1}q+5^{2k+1}q^3-7^{2k+1}q^6+\cdots} {1-3q+5q^3-7q^6+\cdots}
=\frac{\sum_{n\geq0}(-1)^n(2n+1)^{2k+1}q^{\frac{n(n+1)}2}}{\sum_{n\geq0}(-1)^n(2n+1)q^{\frac{n(n+1)}2}}.
 \end{equation}
In terms of the renormalized Eisenstein series
\begin{equation}
E_{2j}(\tau):=\frac{2j}{B_{2j}}\cdot G_{2j}(\tau) = 1 -\frac{4j}{B_{2j}} 
\sum_{n=1}^{\infty}\sigma_{2j-1}(n)q^n,
\end{equation}
Ramanujan found that
$$
U_0=1,\ \ U_2=E_2,\ \ U_4=\frac13(5E_2^2-2E_4),\ \ U_6=\frac19(35E_2^3-42E_2E_4+16E_6),\dots
$$
and he conjectured that every $U_{2k}$ has such an expression. Two of the authors and Singh proved  (see Theorem~1.2 of \cite{AOS2})  this claim, and offered formulas as traces of partition Eisenstein series.

To relate the $\mathcal{\widehat{A}}_k(\tau)$ to Ramanujan's $U_{2k},$ viewed as $q$-series,  we do not use the expressions in Theorem~\ref{QuasimodularGenera}~(1). Instead, we use $E$-normalized traces of partition Eisenstein series
\begin{equation}
\Tr_k^{(E)}(\phi;\tau):=\sum_{\lambda \vdash k} \phi(\lambda)E_{\lambda}(\tau),
\end{equation}
where $E_{\lambda}$ is defined as in  (\ref{primary}), with the Eisenstein series $E_{2j}$ replacing the $G_{2j}.$

It turns out that the  quasimodular $\mathcal{\widehat{A}}_k(\tau)$  are ``partition twists'' of the  $E$-traces of the function\footnote{For aesthetics, we slightly alter the function $\phi_U$ from \cite{AOS2}.}
\begin{eqnarray}
\phi_U(\lambda):=\prod_{j=1}^k\frac1{m_j!}\left(\frac{B_{2j}}{(2j)(2j)!}\right)^{m_j},
\end{eqnarray}
that give Ramanujan's $U_{2k}$ series.

\begin{theorem}\label{RamanujanTwist} 
If $k$ is a positive integer, then as Fourier series the following are true.

(1) We have that
$$
U_{2k}(q)= 4^k (2k+1)! \cdot \Tr_k^{(E)}(\phi_U; \tau).
$$

(2) We have that
$$
\mathcal{\widehat{A}}_k(\tau)= (-1)^k\cdot  \Tr_k^{(E)}(|\phi_U|; \tau).
$$
\end{theorem}

\begin{tworemarks} \ \ \newline
(1) As polynomials  in $E_{\lambda},$  the signs in $\mathcal{\widehat{A}}_k(\tau)$ are the same and are given by $(-1)^k.$
 \newline
(2) Theorem~\ref{RamanujanTwist} shows that Ramanujan's $U_{2k}(q)$ and the $\mathcal{\widehat{A}}_k(\tau)$-genus agree up to choices of sign in the  monomials and explicit scalar multiplier. In particular, the signs differ precisely for those monomials that correspond to $\lambda\vdash k$ with an odd number of parts. 
\end{tworemarks}

\begin{example} Ramanujan's $U_6$ and the $\mathcal{\widehat{A}}_3$-genus are
\begin{displaymath}
\begin{split}
U_6(q)= \frac{16E_6-42E_2E_4+35E_2^3}9  \ \ \ {\text {\rm and}}\ \ \ 
\mathcal{\widehat{A}}_3(\tau)&=\frac{-16E_6-42E_2E_4-35E_2^3}{2903040}. 
\end{split}
\end{displaymath}
The signs differ for the monomials $E_6$ and $E_2^3,$ which correspond to the partitions $\lambda=(3)$ and $\lambda=(1,1,1),$ the partitions of 3 with an odd number of parts. 

 Here we offer a few more examples
 \begin{displaymath}
\begin{split}
&\widehat{\mathcal{A}}_1(\tau)=-\frac{E_2}{24}, \ \ \ \ \ \  \ \ 
\widehat{\mathcal{A}}_2(\tau)=\frac{2E_4+5E_2^2}{5760},\ \ \ \  \ \  \ \ 
\widehat{\mathcal{A}}_3(\tau)=\frac{-16E_6-42E_2E_4-35E_2^3}{2903040},\\ \ \ \\
&\widehat{\mathcal{A}}_4(\tau)=\frac{144E_8+320E_2E_6+84E_4^2+420E_2^2E_4+175E_2^4}{1393459200},\\ \ \ \\
&\widehat{\mathcal{A}}_5(\tau)=\frac{-768E_{10}-1584E_2E_8-704E_4E_6-1760E_2^2E_6-924E_2E_4^2-1540E_2^3E_4-385E_2^5}{367873228800}.
\end{split}
\end{displaymath}
\end{example}

This paper is organized as follows. In Section~2 we prove Theorems~\ref{Theorem1A} and \ref{Theorem1L}, and Corollary~\ref{WittenCorollary} by making use of Weierstrass' theory of elliptic functions and Jacobi forms. In Section~3, we prove Theorem~\ref{QuasimodularGenera} using P\'olya's identity for cycle index polynomials for the symmetric group. We also prove Theorem~\ref{RamanujanTwist} by combining these results with the earlier results from \cite{AOS2}.

\section*{Acknowledgements}
\noindent  The authors thank the two referees for their comments and suggestions. They also thank Jack Morava for helpful conversations and encouragement. They also thank  John Shareshian and Richard Stanley for bringing our attention to Hitchin's article \cite{Hitchin}, and they thank Wei-Lun Tsai for spotting typographical errors in an earlier version of this paper. The third author thanks the Thomas Jefferson Fund and the NSF
(DMS-2002265 and DMS-2055118).

\section{Proof of Theorems~\ref{Theorem1A} and \ref{Theorem1L}}

Here we prove Theorems~\ref{Theorem1A} and \ref{Theorem1L} using the theory of elliptic functions and Jacobi forms. In the next subsection we recall the nuts and bolts that we require about these functions.

\subsection{Jacobi forms and elliptic functions}

We first  recall the definition of a Jacobi form. 
\begin{definition}A holomorphic function $F(z;\tau)$ on $\mathbb{C}\times \mathbb{H}$ is a  {\it Jacobi form for $\SL_2(\Z)$ of weight $k$ and index $m$} if it 
satisfies the following conditions:

\noindent
(1) For all $\gamma=\left(\begin{smallmatrix}a&b\\c&d\end{smallmatrix}\right)\in \SL_2(\Z),$ we have the modular transformation
\[F\left(\frac{z}{c\tau+d};\frac{a\tau+b}{c\tau+d}\right) \ = \ (c\tau+d)^k\exp\left(2\pi i\cdot \frac{mcz^2}{c\tau+d}\right) F(z;\tau).\]

\noindent
(2)  For all integers $a,b$, we have the elliptic transformation
\[
F(z+a\tau+b;\tau)  \ = \  \exp\big(-2\pi i m(a^2 \tau+2az)\big)F(z;\tau).
\]

\noindent
(3) The 
Fourier expansion of $F(z;\tau)$ is given by
\[F(z;\tau)=\sum_{n\geq 0}\sum_{r^2\leq 4mn} b(n,r)q^nu^r,\]
where $b(n,r)$ are complex numbers and $u:=e^{2\pi i z}$.
\end{definition}

As stated in the introduction, the theta function 
\[\theta(z;\tau) = \sum_{n\in \Z} u^nq^{n^2},
\]
 where $u:=e^{2\pi i z}$ and $q:=e^{2\pi i \tau}$ is a Jacobi form of weight $1/2$ and index $1/2$. 
For our purposes, we require and then modify the function
\begin{align}\label{ThetaProdDef}
    \Theta(z; \tau) \ := \  (u^{1/2}-u^{-1/2})\prod_{n\geq 1}\frac{(1-uq^n)(1-u^{-1}q^n)}{(1-q^n)^2}.
\end{align}
This function is related to both the function $\widetilde \Theta(z;\tau)$ defined in (\ref{ThetaDef}), and $\theta(z;\tau),$ as shown below.

\begin{proposition}\label{ThetaFormulas} The following identities are true.

\noindent
(1) In terms of $\theta(z;\tau)$ and Dededkind's eta function $\eta(\tau):=q^{\frac{1}{24}}\prod_{n=1}^{\infty}(1-q^n),$ we have that
\[
\Theta(z;\tau)=\frac{1}{\eta(\tau)^3}\cdot u^{\frac{1}{2}}q^{\frac{1}{8}}\cdot\theta\left (z+\frac{\tau}{2}+\frac{1}{2}; \tau \right).\]
 
 \noindent
 (2) We have that
 \[
\widetilde\Theta(z;\tau)=\exp\left(\frac{\pi}{2}\cdot \frac{z^2}{\Im(\tau)}\right)\Theta(z;\tau).
\]
\end{proposition}
\begin{remark}
Combining the modular transformation properties of Dedekind's eta-function $\eta(\tau)$ (for example, see Chapter~1.4 of \cite{CBMS}) with Proposition~\ref{ThetaFormulas} (1), we have  that  $\Theta(z;\tau)$ is a Jacobi form of weight $-1$ and index $1/2$. 
\end{remark}

\medskip
\noindent
{\it Proof of Proposition~\ref{ThetaFormulas}.} 
Claim (1) follows as an easy application of the Jacobi Triple Product formula (see Theorem 2.8 of \cite{Andrews}), which allows us to write
\[\begin{split}
(u^{1/2}-u^{-1/2})\prod_{n\geq 1}\frac{(1-uq^n)(1-u^{-1}q^n)}{(1-q^n)^2} &= \frac{1}{\eta(\tau)^3}\sum_{n\in \mathbb{Z}} (-1)^n u^{\frac{2n+1}{2}}q^{\frac{(2n+1)^2}{2}}\\
&=\frac{1}{\eta(\tau)^3}\cdot u^{1/2}q^{1/8}\theta(z+\tfrac{\tau}{2}+\tfrac{1}{2};\tau).
\end{split}
\]

The second claim follows immediately from (1) and (\ref{ThetaDef}). $\ \ \square$
\bigskip

To prove Theorems~\ref{Theorem1A} and \ref{Theorem1L}, we require
 the Weierstrass $\sigma$-function,
\begin{eqnarray}\label{sigmafunction}
\sigma(z,\tau):=z\prod_{\substack{w\in \Lambda_\tau\\ \Im (w)>0 \text{ or }w>0 }} \left(1-\frac{z^2}{w^2}\right)\exp\left(\tfrac{z^2}{w^2}\right),
\end{eqnarray}
where $\Lambda_\tau$ is the lattice $\Lambda_\tau= \Z\tau+\Z.$  We have the following elementary proposition.

\begin{proposition}\label{sigma Theta}
We have that
$$
\sigma(z;\tau) = \frac{1}{2\pi i} e^{\frac{G_2(\tau)}{2} (2\pi i z)^2}\cdot \Theta(z;\tau).
$$
\end{proposition}
\medskip

\noindent
{\it Proof.}
The $\sigma$-function has a   $q$-series expansion (see Theorem I.6.6.4 of \cite{Silverman}) given by
\[
\sigma(z;\tau)= \frac{1}{2\pi i} e^{\frac{G_2(\tau)}{2} (2\pi i z)^2}(u^{1/2}-u^{-1/2})\prod_{n\geq 1}\frac{(1-uq^n)(1-u^{-1}q^n)}{(1-q^n)^2}.
\]
Thus $\sigma$ is also related to the modified theta function $\Theta(z;\tau)$, defined by (\ref{ThetaProdDef}), as claimed. $\ \square$
\bigskip

Finally, we will need a lemma giving a slightly nonstandard formula for the weight 2 nonholomorphic  Eisenstein series  $G_2^{\star}(\tau):=1/4\pi\Im(\tau)+G_2(\tau)$

\begin{lemma}\label{G2limit}
We have that 
\[
 (2\pi i)^2G_{2}^\star(\tau)=
\lim_{s\to 0^+} \, \sum_{k=1}^\infty \sum_{\substack{m,n\in \Z\\ \gcd (m,n)=1}}\frac{1}{k^{2}\cdot (m\tau+n)^2|m\tau+n|^s}.
\]
\end{lemma}

\bigskip
\noindent
{\it Proof.}  The standard application of ``Hecke's trick''  (for example, see p. 84 of \cite{BFOR}), to force convergence of the weight $2$ Eisenstein series, gives the formula
\[
 (2\pi i)^2G_{2}^\star(\tau)=\lim_{s\to 0^+}\sum_{\substack{m,n\in \Z\\ (m,n)\neq (0,0)}}\frac{1}{(m\tau+n)^2|m\tau+n|^s}.
\]
The expression in the limit factors as 
\[
\begin{split}
\sum_{\substack{m,n\in \Z\\ (m,n)\neq (0,0)}}\frac{1}{(m\tau+n)^2|m\tau+n|^s} &=\sum_{k=1}^\infty \sum_{\substack{m,n\in \Z\\ \gcd(m,n)=1}}\frac{1}{k^{2+s}(m\tau+n)^2|m\tau+n|^s}\\
& = \zeta(2+s)\cdot \sum_{\substack{m,n\in \Z\\ \gcd(m,n)=1}}\frac{1}{(m\tau+n)^2|m\tau+n|^s},
\end{split}
\]
where in each term we have factored out $k=\gcd(m,n)$, and $\zeta(s)$ is the Riemann zeta-function. Similarly, we have that 
\[
\sum_{k=1}^\infty\sum_{\substack{m,n\in \Z\\ (m,n)\neq (0,0)}}\frac{1}{k^2\cdot (m\tau+n)^2|m\tau+n|^s}  = \zeta(2)\cdot \sum_{\substack{m,n\in \Z\\ \gcd(m,n)=1}}\frac{1}{(m\tau+n)^2|m\tau+n|^s}.
\]
The lemma follows since $\zeta(s)$ is continuous at $2$. $\square$
\medskip

\subsection{Proof of Theorem~\ref{Theorem1A}}

We first find the Weierstrass factorization of the characteristic series (see (\ref{Gen_sinc}))
$$Q_{\widehat A}(x)=\frac{\sqrt{x}}{\sinh(\sqrt{x})}.
$$
The function $\sin(x)/x$ has the well-known Weierstrass factorization 
\[
\frac{\sin(x)}{x} = \prod_{k=1}^\infty\left(1-\frac{x^2}{\pi^2 k^2}\right).
\]
This gives the factorization for $Q_{\widehat A}(z)$ by applying the identity $\sinh(x)=\sin(-ix).$ 

Using this factorization, we have that
\[
\begin{split}
\widehat A (X_\tau(s);(2\pi i z)^2) &= \prod_{k=1}^\infty\prod_{x\in X_\tau(s)}\left(1+\frac{(2\pi i z)^2 x}{4 \pi^2 k^2}\right)^{-1}\\
&= \prod_{k=1}^\infty\prod_{x\in X_\tau(s)}\left(1-\frac{ z^2 \cdot x}{k^2}\right)^{-1}\\
&=\prod_{k=1}^\infty\prod_{x\in X_\tau(s)}\left(1-\frac{ z^2 \cdot x}{k^2}\right)^{-1}\exp\left( -\frac{z^2\cdot x}{k^2}+\frac{z^2\cdot x}{k^2}\right).
\end{split}
\]
The last step allows us to break the expression in two parts, which behave differently as $s\to 0^+$.

 For the first piece, we may simply evaluate at $s=0$ and use Proposition~\ref{sigmafunction} to obtain
 \[
 \begin{split}
\lim_{s\to 0^+}\prod_{k=1}^\infty\prod_{x\in X_\tau(s)}&\left(1-\frac{ z^2 \cdot x}{k^2}\right)^{-1}\exp\left( -\frac{z^2\cdot x}{k^2}\right)\\ &= \prod_{k=1}^\infty\prod_{\substack {(m,n) \in \Z^2/(\pm 1)\\ \gcd(m,n)=1}}\left(1-\frac{ z^2 }{k^2(m\tau+n)^2}\right)^{-1}\exp\left( -\frac{z^2}{k^2\cdot(m\tau+n)^2}\right)\\
&=\frac{z}{\sigma (z;\tau)}.
\end{split}
\]
For the second piece, we use Lemma \ref{G2limit} to obtain
\[\begin{split}
\lim_{s\to 0^+}\exp\left(\frac{z^2\cdot x}{k^2}\right)&=  \lim_{s\to 0^+}
\prod_{k=1}^\infty \prod_{\substack {(m,n) \in \Z^2/(\pm 1)\\ \gcd(m,n)=1}}
\exp\left( \frac{z^2}{k^2\cdot (m\tau+n)^2|m\tau+n|^s}\right)\\
& = \exp\left( \frac{1}{2}G_2^\star(\tau) (2\pi i z)^2\right).
\end{split}
\]
Here the $1/2$ appears since this expression is a sum over only a half-lattice, whereas  this lemma uses the sum over the full lattice. 

Using Proposition~\ref{sigma Theta} and Proposition \ref{ThetaFormulas}, we obtain the claimed expression 
\[
\lim _{s\to 0} \widehat A (X_\tau(s);(2\pi i z)^2) = \frac{2\pi i z}{\widetilde \Theta(z;\tau)}.
\]

\subsection{Proof of Theorem~\ref{Theorem1L}}
Following the proof of Theorem~\ref{Theorem1A}, we first find the Weierstraas factorization of the characteristic series (see (\ref{LGenus}))
$$\Q_L(x)=\frac{\sqrt{z}}{\tanh(\sqrt{z})}.
$$
The function $\tan(x)/x$ has Weierstrass factorization 
\[
\frac{\tan(x)}{x} = \frac{\prod_{k=1}^\infty\left(1-\frac{x^2}{\pi^2 k^2}\right)}{\prod_{k=1}^\infty\left(1-\frac{4 x^2}{\pi^2 (2k-1)^2}\right)}.
\]

Using (\ref{LGenus}), and the fact that $\tanh(x)=\tan(-ix),$ we have that 

\[\begin{split}
L(X_\tau(s),(\pi i z)^2) &= \prod_{x\in X_\tau(s)}\frac{\prod_{k=1}^\infty\left(1+\frac{4 (\pi i z)^2 x }{\pi^2 (2k-1)^2}\right)}{\prod_{k=1}^\infty\left(1+\frac{(\pi i z
)^2 x}{\pi^2 k^2}\right)}\\
&=\prod_{x\in X_\tau(s)}\frac{\prod_{k=1}^\infty\left(1-\frac{4 z^2 x }{(2k-1)^2}\right)}{\prod_{k=1}^\infty\left(1-\frac{ z^2 x}{k^2}\right)} 
=\prod_{x\in X_\tau(s)}\frac{\prod_{k=1}^\infty\left(1-\frac{4 z^2 x }{k^2}\right)}{\prod_{k=1}^\infty\left(1-\frac{ z^2 x}{k^2}\right)^2}.
\end{split}
\]
Using the calculations from the previous subsection again, 
we find that 
\[
\lim _{s\to 0^+}L(X_\tau(s),(\pi i z)^2) = \pi i z\cdot \frac{\widetilde \Theta(2z;\tau)}{\widetilde \Theta(z;\tau)^2}.
\]

\subsection{Proof of Corollary~\ref{WittenCorollary}}
The logarithmic derivative of the $\sigma$-function (with respect to $z$) has Taylor expansion
\[
\frac{\sigma'(z;\tau)}{\sigma(z;\tau)} = \frac{1}{z} -\sum_{k\geq2} \frac{G_{2k}(\tau)(2\pi i)^{2k}}{(2k-1)!}z^{2k-1}
\]
(see Prop. I.5.1 of \cite{Silverman}, where we note a difference in notation with our $G_{2k}(\tau)$ being $\frac{(2k-1)!}{(2\pi i)^{2k}}G_{2k}(\Lambda_\tau)$).
This gives us the exponential expansion of $\sigma$ as
 \[
 \sigma(z;\tau) = z\cdot \exp\left( -\sum_{k\geq 2}\frac{G_{2k}(\tau)}{(2k)!}(2\pi i z)^{2k}\right),
 \] 
 Therefore, we find that the characteristic series of the Witten genus is (see (\ref{WittenGenus})) satisfies
\[
Q_W(z)=\exp\left(\sum_{k\geq 2}\frac{G_{2k}(\tau)(2\pi i z)^{2k}}{(2k)!}\right) = \frac{z}{\sigma(z;\tau)}.
\]
Combining Theorem \ref{Theorem1A}, Proposition~\ref{ThetaFormulas}, and Proposition~\ref{sigma Theta} we obtain the claim.

\section{Proof of Theorems~\ref{QuasimodularGenera} and \ref{RamanujanTwist}}

In this section we prove Theorems~\ref{QuasimodularGenera} and \ref{RamanujanTwist}, which express the quasimodular representations of the $\widehat{A}$-genus and $L$-genus as traces of partition Eisenstein series. To obtain these results, we make use of exponential generating functions that arose in the previous section. These generating functions can be reformulated as traces of partition Eisenstein series using special identities within the framework of P\`olya's theory of cycle index polynomials.

\subsection{P\`olya's cycle index polynomials}

The structure of traces of partition Eisenstein series arises from the classical theory of the symmetric group, and their connection to integer partitions. Namely, the key tool is P\'olya's theory of cycle index polynomials  (for example, see \cite{Stanley}). Recall that  a partition $\lambda=(\lambda_1,\dots,\lambda_{\ell(\lambda)})\vdash k$ or $(1^{m_1},\dots,k^{m_k})\vdash k$, labels a conjugacy class by cycle type. Moreover, the number of permutations in $\mathfrak{S}_k$ of cycle type $\lambda$ is
is $k!/z_{\lambda}$, where 
 $z_{\lambda}:=1^{m_1}\cdots k^{m_k}m_1!\cdots m_k!$. The {\it cycle index polynomial} for the symmetric group $\mathfrak{S}_k$ is given by
\begin{align} \label{CIF}
Z(\mathfrak{S}_k)&=\sum_{\lambda\vdash k}\frac1{z_{\lambda}}\prod_{j=1}^{\ell(\lambda)}x_{\lambda_j}
=\sum_{\lambda\vdash k}\prod_{j=1}^k\frac1{m_j!}\left(\frac{x_j}{j}\right)^{m_j}.
\end{align}
We require the following generating function for these polynomials in $k$-aspect.
\begin{lemma}[Example 5.2.10 of \cite{Stanley}]\label{PolyaGenFunction} As a power series in $y$, the generating function for the cycle index polynomials satisfies
$$\sum_{k\geq0}Z(\mathfrak{S}_k)\,t^k=\exp\left(\sum_{k\geq1} x_k\cdot \frac{t^k}k \right).$$
\end{lemma}

\begin{example}
 Here are the first few examples of P\'olya's cycle index polynomials: 
$$Z(\mathfrak{S}_1)=x_1, \qquad Z(\mathfrak{S}_2)=\frac1{2!}(x_1^2+x_2), \qquad Z(\mathfrak{S}_3)=\frac1{3!}(x_1^3+3x_1x_2+2x_3).
$$
\end{example}

\subsection{Proof of Theorem~\ref{QuasimodularGenera}}
The characteristic series $Q_{\widehat{A}}(z)$ (see \eqref{Gen_sinc}) has a well-known (see \cite[eq. (3.1)]{AOS2}) explicit exponential generating function
\begin{equation}\label{QaNice}
Q_{\widehat{A}}(z)=\frac{\sqrt{z}/2}{\sinh(\sqrt{z}/2)}=\exp\left(-\sum_{j=1}^{\infty} \frac{B_{2j}z^j}{(2j)(2j)!}\right),
\end{equation}
which enables us, by (\ref{infinite_product}), to obtain
$$\widehat{A}(s_1,s_2,\dots;t)=\prod_{i=1}^{\infty}Q_{\widehat{A}}(x_it)=\exp\left(\sum_{k=1}^{\infty} \frac{-B_{2k}s_kt^k}{(2k)(2k)!}\right).$$
Here we see the natural role of the  power sum symmetric functions $\{s_k\}$. 

To prove the theorem, we invoke P\'olya's formula in Lemma~\ref{PolyaGenFunction}, with $\frac{w_k}k=\frac{-B_{2k}s_k}{(2k)(2k)!}.$  In this way, we obtain
\begin{displaymath}
\begin{split}
\widehat{A}(s_1,s_2,\dots;t)
&=\sum_{k\geq0} t^k\sum_{\lambda\vdash k} \prod_{j=1}^k \frac1{m_j!} \left(\frac{-B_{2j}s_j}{(2j)(2j)!}\right)^{m_j}
\\
&=\sum_{k\geq0} t^k \sum_{\lambda\vdash k} s_{\lambda} \prod_{j=1}^k \left(\frac{B_{2j}}{2j}\right)^{m_j} \cdot
 \prod_{j=1}^k \frac1{m_j!} \left(\frac{-1}{(2j)!}\right)^{m_j}.
\end{split}
\end{displaymath}
Under the identification $B_{2j} s_j/2j  \longleftrightarrow G_{2j}$, we have the desired expression
as a trace of partition Eisenstein series. Namely, we find that  $$\widehat{\mathcal{A}}_k(\tau)= \Tr_k(\phi_{\widehat{A}};\tau).$$

Now we turn to the case of the $L$-genus, which has characteristic series (see (\ref{LGenus}))
$$
Q_L(z)=\frac{\sqrt{z}}{\tanh(\sqrt{z})}.
$$
On the other hand, one may recall the series exapansion (easily derived from that of $\tan x$)
$$\cosh(\sqrt{z})=\exp\left(\sum_{j=1}^{\infty} \frac{4^j(4^j-1)B_{2j}z^j}{(2j)(2j)!}\right).$$
Combining this with formula (\ref{QaNice}) for $Q_{\widehat{A}}(z)$, we get
$$Q_L(z)=\frac{\sqrt{z}}{\tanh(\sqrt{z})}=\frac{\sqrt{z}}{\sinh(\sqrt{z})}\cdot \cosh(\sqrt{z})
=\exp\left(\sum_{j=1}^{\infty} \frac{4^j(4^j-2)B_{2j}z^j}{(2j)(2j)!}\right).$$
Arguing as above with (\ref{infinite_product}) and P\'olya's Lemma~\ref{PolyaGenFunction} {\it mutatis mutandis}, we obtain the claimed conclusion
$$
\mathcal{L}_k(\tau)= \Tr_k(\phi_L;\tau).
$$

\subsection{Proof of Theorem~\ref{RamanujanTwist}}

Claim (1) is a simple reformulation of Theorem~1.2 (1) of \cite{AOS2}. The reader merely needs to be aware of the different normalizations of the function $\phi_U$.

The proof of claim (2) is a little more involved.
Beginning with \eqref{phi_A} and Theorem~\ref{QuasimodularGenera}, we apply the correspondence $G_{2j} \longleftrightarrow B_{2j}s_j/2j$ as follows
\begin{align*} \mathcal{\widehat{A}}_k(\tau) &= \sum_{\lambda\vdash k} G_{\lambda}(\tau) \prod_{j=1}^k\frac1{m_j!}\left(\frac{-1}{(2j)!}\right)^{m_j}
=\sum_{\lambda\vdash k} E_{\lambda}(\tau) \prod_{j=1}^k\frac1{m_j!}\left(\frac{-B_{2j}}{(2j)(2j)!}\right)^{m_j}.
\end{align*}
Since $B_{2j}=(-1)^{j-1}\vert B_{2j}\vert$ (or $-B_{2j}=(-1)^j\vert B_{2j}\vert$) and $\sum_{j=1}^kjm_j=k$, it follows that
$$\mathcal{\widehat{A}}_k(\tau)
=\sum_{\lambda\vdash k} (-1)^{\sum_j jm_j} E_{\lambda} \prod_{j=1}^k\frac1{m_j!}\left(\frac{\vert B_{2j}\vert}{(2j)(2j)!}\right)^{m_j}
=(-1)^k \sum_{\lambda\vdash k} E_{\lambda}\cdot \vert \phi_U\vert=(-1)^k\Tr_k^{(E)}(\vert\phi_U\vert;\tau).$$
This proves the desired expression in claim (2) as a twisted $E$-trace of partition Eisenstein series.


\begin{thebibliography}{99}

\bibitem{Hirzebruch} F. Hirzebruch, \emph{Topological methods in algebraic geometry}, Classics in Mathematics, Springer-Verlag, Berlin, 1995.

\bibitem{HBJ} F. Hirzebruch, T. Berger, and R. Jung, \emph{Manifolds and modular forms}, Aspects of Math. E20, Friedr. Vieweg \& Sohn, Braunschweig, 1992.

\bibitem{BorelHirzebruch} A. Borel and F. Hirzebruch, \emph{Characteristic classes and homogeneous spaces. II}, Amer. J. Math.
\textbf{81} (1959), 315--382.

\bibitem{Hitchin} N. Hitchin, \emph{The Dirac operator}, Bull. Amer. Math. Soc.  \textbf{67} no. 1 (2025), 3-16.


\bibitem{BFOR} K. Bringmann, A. Folsom, K. Ono, and L. Rolen, \emph{Harmonic Maass forms and modular forms: theory and applications}, Amer. Math. Soc. Colloq. Series \textbf{64}, Amer. Math. Soc., Providence, 2017.


\bibitem{Zagier} M. ~Kaneko and D. ~Zagier,
\emph{A generalized Jacobi theta function and quasimodular forms},
The moduli space of curves (Texas Island, 1994), Progr. Math., vol. \textbf{129}, Birkh\"auser Boston, Boston, MA, 1995, 165--172.

\bibitem{EichlerZagier} M. Eichler and D. Zagier, \emph{The theory of Jacobi forms}, Prog. Math. \textbf{55}, Birkh\"auser, Boston, 1985.

\bibitem{AGOS} T. Amdeberhan, M. Griffin, K. Ono, and A. Singh, \emph{Traces of partition Eisenstein series}, Forum Mathematicum, \textbf{37} (2025), no. 6, 1835--1847.

\bibitem{CBMS}  K. Ono,
 \emph{The web of modularity: arithmetic of the coefficients of modular forms and $q$-series}, CBMS Regional Conference Series in Mathematics, \textbf{102}, Amer. Math. Soc., Providence, RI, 2004.

\bibitem{Witten} E. Witten, \emph{The index of the Dirac operator in loop space}, Elliptic curves and modular forms in Algebraic Topology, Princeton 1986, Springer Lect. Notes Math. \textbf{1326}, 1988, 161-181.


\bibitem{Hopkins} M. J. Hopkins, \emph{Topological modular forms, the Witten genus, and theorem of the cube}, Proc. Intl. Congress Math. \textbf{1} (Z\"urich, 1994), Birkh\"auser, Basel, 1995, 554--565.


\bibitem{AOS} T. Amdeberhan, K. Ono, and A. Singh,
\emph{MacMahon's sums-of-divisors and allied q-series}, Adv. Math. \textbf{452} (2024), Article 109820.

\bibitem{AOS2} T. Amdeberhan, K. Ono, and A. Singh, \emph{Derivatives of theta functions as traces of partition Eisenstein series},
Nagoya Mathematical Journal, (2025): 1--12 (doi:10.1017/nmj.2024.30).


\bibitem{Andrews} G. E. Andrews, \emph{Theory of partitions}, Cambridge Univ. Press, 1998.

\bibitem{JustSchneider} M. Just and R. Schneider, \emph{Partition Eisenstein series and semi-modular forms}, Res. Number Theory {\textbf 7} (2021), Paper No. 61.

\bibitem{Rama} S. ~Ramanujan, 
\emph{The lost notebook and other unpublished papers}, (1988)
New Delhi; Berlin, New York: Narosa Publishing House; Springer-Verlag, Reprinted (2008).


\bibitem{Silverman} J. H.~Silverman,
\emph{Advanced topics in the arithmetic of elliptic curves}, 
Grad. Texts in Math., \textbf{151} Springer-Verlag, New York (1994) xiv+525 pp.


\bibitem{Stanley} R. P. Stanley, \emph{Enumerative combinatorics. Vol. 2}, Cambridge Studies in Math. {\bf 62} Cambridge Univ. Press, 1999.


\end{thebibliography}
\end{document}